\documentclass[11pt]{amsart}
\usepackage{amsmath,amsfonts,amstext,amsthm,epsfig}

\newtheorem{thm}{Theorem}

\theoremstyle{remark}
\newtheorem{rem}{Remark}
\numberwithin{equation}{section}

\newcommand{\eps}{\varepsilon}

\def\R{\mathbb R}

\begin{document}

\title[A remark on well-posedness...]
{A remark on well-posedness for hyperbolic equations with singular
coefficients}%
\author{Daniele Del Santo}\author{Martino Prizzi}

\address{Daniele Del Santo, Universit\`a di Trieste, Dipartimento di
Scienze Matematiche,
Via Valerio 12/1, 34127 Trieste, Italy}%
\email{delsanto@univ.trieste.it}%

\address{Martino Prizzi, Universit\`a di Trieste, Dipartimento di
Scienze Matematiche,
Via Valerio 12/1, 34127 Trieste, Italy}%
\email{prizzi@mathsun1.univ.trieste.it}%
\thanks{We would like to thank Prof.~F. Colombini for the stimulating and
fruitful conversations.}%

\begin{abstract}
We prove some $C^\infty$ and Gevrey well-posedness results for hyperbolic
equations with singular coefficients.
\end{abstract}
\maketitle

\section{Introduction}

This work is devoted to the study of the well-posedness of the Cauchy
problem for a linear hyperbolic operator whose coefficients depend only on
time.

We consider the equation
\begin{equation}\label{eq1}
u_{tt}-\sum_{i,j=1}^n a_{ij}(t)u_{x_ix_j}=0
\end{equation}
in $[0,T]\times\R^n$, with initial data
\begin{equation}\label{eq2}
u(0,x)=u_0(x),\quad u_t(0,x)=u_1(x)
\end{equation}
in $\R^n$. The matrix $(a_{ij})$ is supposed to be real and symmetric. Setting
\begin{equation}\label{eq3}
a(t,\xi):=\sum_{i,j=1}^n a_{ij}(t)\xi_i\xi_j/|\xi|^2,\quad(t,\xi)\in[0,T]
\times(\R^n\setminus \{0\}),
\end{equation}
we assume that $a(\cdot,\xi)\in L^1(0,T)$ for all $\xi\in\R^n\setminus\{0\}$.

We suppose that the equation (\ref{eq1}) is hyperbolic i.e.
\begin{equation}\label{eq4}
a(t,\xi)\geq \lambda_0\geq 0
\end{equation}
for all $(t,\xi)\in[0,T]\times(\R^n\setminus \{0\})$.

In the strictly hyperbolic case (i.e. $\lambda_0>0$) it is well known that
if the
coefficients $a_{ij}$ are Lipschitz-continuous then the Cauchy problem
(\ref{eq1}),
(\ref{eq2}) is well-posed in Sobolev spaces. In the same case if the $a_{ij}$'s
are Log-Lipschitz-continuous or H\"older-continuous of index $\alpha$,
(\ref{eq1}),
(\ref{eq2}) is well-posed in $C^\infty$ or in the Gevrey space $\gamma^{(s)}$ for 
$s<{1\over 1-\alpha}$ respectively (see \cite{CDGS}). 
In the weakly hyperbolic case (i.e. $\lambda_0=0$) if the coefficients are
$C^{k,\alpha}$ then the problem (\ref{eq1}), (\ref{eq2}) is
$\gamma^{(s)}$-well-posed
for $s<1+{k+\alpha\over 2}$ (see \cite{CJS}). Some counter examples show
that all these results are sharp (see also \cite{CL}).

Recently Colombini, Del Santo and Kinoshita have considered the same
problem for
operators having coefficients which are $C^1$ on  $[0,T]\setminus\{t_0\}$
with a
singularity concentrated at $t_0$. In this situation, under the main
assumptions that
\begin{equation}
\begin{array}{l}
|t_0-\cdot|^pa'(\cdot,\xi)=\beta(\cdot,\xi)
\in L^\infty (0,T)\quad {\rm  for\ all}\ \xi\in\R^n\setminus\{0\}\\[0.3cm]
|t_0-\cdot|^ra(\cdot,\xi)=\alpha(\cdot,\xi)
\in L^\infty (0,T)\quad {\rm  for\ all}\ \xi\in\R^n\setminus\{0\}\\
\end{array}
\label{eq5}
\end{equation}
it is possible to show that the Cauchy problem (\ref{eq1}), (\ref{eq2}) is
$\gamma^{(s)}$-well-posed, the value of $s$ depending on $p$ and $r$ (see
\cite{CDSK1}
and\cite{CDSK2}) (here and in the following ``$\,'\,$" denotes the
differentiation
with respect to $t$).

The aim of the present work is to improve the results of \cite{CDSK1} and
\cite{CDSK2} allowing the function $\beta$ in (\ref{eq5}) to be in a $L^q$
space and
removing the growth assumption on $a$. We make the following assumptions: let
$1\leq q\leq+\infty$ and
$p\geq 0$ and let
$t_0\in[0,T]$;  suppose that
\vskip 0.3cm
\begin{enumerate}
\item[(H1)] $a(\cdot,\xi)\in
\cap_{\eps>0}W^{1,1}(]0,t_0-\eps[\cup]t_0+\eps,T[)$
for all $\xi\in\R^n\setminus\{0\}$;\vskip 0.3cm
\item[(H2)] $|t_0-\cdot|^pa'(\cdot,\xi)=\beta(\cdot,\xi)
\in L^q(0,T)$ for all $\xi\in\R^n\setminus\{0\}$.
\end{enumerate}
\vskip 0.3cm

In the weakly hyperbolic case the results are the following.
\begin{thm}\label{th1}
Assume that $3\leq (p+1/q)$.
Then the Cauchy problem
(\ref{eq1}), (\ref{eq2}) is $\gamma^{(\sigma)}$-well-posed for
$1\leq\sigma<
\frac{(p+1/q)-\frac32}{(p+1/q)-2}$. If moreover
\begin{equation}
|t_0-\cdot|^ra(\cdot,\xi)=\alpha(\cdot,\xi)
\in L^s(0,T)\quad\text{ for all $\xi\in\R^n\setminus\{0\}$,}
\label{supplcond}
\end{equation}
with $r\geq 0$, $1\leq s\leq+\infty$ and $(r+1/s)\leq1$, then the Cauchy
problem
(\ref{eq1}), (\ref{eq2}) is $\gamma^{(\sigma)}$-well-posed for
$1\leq\sigma<
\frac{(p+1/q)-\frac32(r+1/s)}{(p+1/q)-(r+1/s)-1}$.
\end{thm}

\begin{thm}\label{th2}
Assume that $(p+1/q)<3$. Then the Cauchy problem
(\ref{eq1}), (\ref{eq2}) is $\gamma^{(\sigma)}$-well-posed for all
$1\leq\sigma<\frac32$.
\end{thm}

The result concerning the strictly hyperbolic case are contained in the
following
theorems.

\begin{thm}\label{th3}
Assume that $1<(p+1/q)<3$. Moreover, assume that $\lambda_0>0$.
Then the Cauchy problem
(\ref{eq1}), (\ref{eq2}) is $\gamma^{(\sigma)}$-well-posed for all
$1\leq\sigma<\frac{(p+1/q)}{(p+1/q)-1}$.
\end{thm}

\begin{thm}\label{th4}
Assume that $(p+1/q)\leq 1$. Moreover, assume that $\lambda_0>0$.
Then the Cauchy problem
(\ref{eq1}), (\ref{eq2}) is $C^\infty$-well-posed.
\end{thm}

\begin{rem}{\rm Adapting to the present situation some counter examples
contained
in \cite{CJS}, \cite{CDSK1}, and \cite{CDSK2} it is possible to see that
the results
of Theorems~\ref{th1}--\ref{th4} are optimal. Let us show this in some
detail in the
case of Theorem \ref{th1}. Suppose $p_0+1/q_0=3$. In this case
$\frac{(p_0+1/q_0)-\frac32}{(p_0+1/q_0)-2}=
\frac{(p_0+1/q_0)-\frac32(r_0+1/s_0)}{(p_0+1/q_0)-(r_0+1/s_0)-1}=\frac32$;
consequently Theorem~2 in \cite{CJS} shows that this value of the Gevrey index
cannot be improved.
Consider next the case that $p_0+1/q_0>3$ and $(r_0+1/s_0)\leq 1$. Let
$\bar\sigma > \sigma_0
=\frac{(p_0+1/q_0)-\frac32(r_0+1/s_0)}{(p_0+1/q_0)-(r_0+1/s_0)-1}$. We fix
$q_1>q_0$
and $s_1>s_0$ in such a way that $p_0+1/q_1>3$, $r_0+1/s_1<1$ and
$\sigma_0<\sigma_1:=\frac{(p_0+1/q_1)-\frac32(r_0+1/s_1)}{(p_0+1/q_1)-(r_0+1/s_1)
-1}
<\bar\sigma$. From Theorem~4 in \cite{CDSK2} we have that there exists a
function
$a:[0,1[\, \to [1/2,+\infty[$ such that $a\in C^\infty([0,1[)$ and
$$
(1-t)^{p_0+1/q_1}a'(t)\in L^\infty, \quad (1-t)^{r_0+1/s_1}a(t)\in L^\infty,
$$
and there exist $u_0$, $u_1\in \gamma^{(\sigma)}$ for all $\sigma>\sigma_1$
such that
the Cauchy problem
\begin{equation}
u_{tt}-a(t)u_{xx}=0, \qquad u(0,x)=u_0(x),\quad u_t(0,x)=u_1(x),
\label{ne}
\end{equation}
has no solution in $W^{2,1}([0,1],{{\mathcal D}'}^{(\sigma)}(\R))$ for all
$\sigma>\sigma_1$.
Consequently
$$
(1-t)^{p_0}a'(t)\in L^{q_0}, \quad (1-t)^{r_0}a(t)\in L^{s_0},
$$
$u_0$, $u_1\in \gamma^{(\bar\sigma)}$ and the Cauchy problem (\ref{ne})
does not have
a solution in $W^{2,1}([0,1],{{\mathcal D}'}^{(\bar\sigma)}(\R))$.}
\end{rem}
\begin{rem}
{\rm Let us remark that Theorem~\ref{th1} is a nontrivial improvement of
Theorem~2 in
\cite{CDSK2} also in the case of $q=\infty$. In fact the growth condition
on $a$ is removed and the result is sharp (see \cite[Th. 4]{CDSK2}).}
\end{rem}

\section{Proof of Theorems \ref{th1}--\ref{th4}}
As a preliminary step, let us observe that, since the coefficients $a_{ij}$
are real integrable functions, the Cauchy problem (\ref{eq1}), (\ref{eq2})
is well posed in
${\mathcal A}'(\R^n)$, the space of real analytic functionals. Moreover, if
the
initial data vanish in a ball, then the solution vanishes in a cone, whose
slope
depends on the coefficients $a_{ij}$.  Therefore it will be sufficient to
show that,
under the hypotheses of each theorem, if $u_0$ and
$u_1$ have compact support then the corresponding solution
$u$ is not only in $W^{2,1}([0,T],{\mathcal A}'(\R^n))$, but it belongs
to a Gevrey space in the $x$ variable. Our main tools in doing this
will be the Paley-Wiener theorem (in the version of
\cite[p. 517]{CDGS}, to which we refer here and throughout) and some energy estimates.

Denoting by $v$ the Fourier transform of $u$ with respect to $x$, equation
(\ref{eq1}) reads
\begin{equation}\label{eq2.1}
v''(t,\xi)+a(t,\xi)|\xi|^2v(t,\xi)=0.
\end{equation}
Let $\epsilon$ be a positive parameter and for each
$\epsilon$ let $a_\eps\colon[0,T]\times(\R^n\setminus\{0\})\to\R$ be a
strictly positive
real function
such that $a_\eps(\cdot,\xi)\in W^{1,1}(0,T)$
for all $\xi\in\R^n\setminus\{0\}$.
We define the {\em approximate energy} of $v$ by
\begin{equation}\label{eq2.2}
E_\eps(t,\xi):=a_\eps(t,\xi)|\xi|^2|v(t,\xi)|^2+|v'(t,\xi)|^2,\quad (t,\xi)\in
[0,T]\times(\R^n\setminus\{0\}).
\end{equation}
Differentiating $E_\eps$ with respect to $t$ and using (\ref{eq2.1}) we get
\begin{multline*}
E'_\eps(t,\xi)=a'_\eps(t,\xi)|\xi|^2|v(t,\xi)|^2+
2a_\eps(t,\xi)|\xi|^2{\rm Re}(v'(t,\xi)\bar v(t,\xi))\\
+2{\rm Re} (v''(t,\xi)\bar v'(t,\xi))\\
\leq \left(\frac{|a'_\eps(t,\xi|}{a_\eps(t,\xi)}+
\frac{|a_\eps(t,\xi)-a(t,\xi)|}{a_\eps(t,\xi)^{1/2}}|\xi|\right)E_\eps(t,\xi).
\end{multline*}
By Gronwall's lemma we obtain
\begin{equation}\label{eq2.3}
E_\eps(t,\xi)\leq E_\eps(0,\xi)
\exp\left(\int_0^T\frac{|a'_\eps(t,\xi|}{a_\eps(t,\xi)}\,d t+
|\xi|\int_0^T\frac{|a_\eps(t,\xi)-a(t,\xi)|}{a_\eps(t,\xi)^{1/2}}\,d t\right)
\end{equation}
for all $t\in[0,T]$ and for all $\xi\in\R^n$, $|\xi|\geq1$.

Now we are able to give the
\begin{proof}[Proof of Theorem \ref{th1}]
First of all, observe that condition (\ref{supplcond}) is always satisfied
at least with $r=0$ and $s=1$ (recall that $a_{ij}\in L^1(0,T)$).

Since $u_0$, $u_1\in\gamma^{(\sigma)}\cap C^\infty_0$, the Paley-Wiener
theorem  ensures that there exist $M,\delta>0$ such that
\begin{equation}\label{eq2.4}
|v(0,\xi)|^2+|v'(0,\xi)|^2\leq M\exp(-\delta|\xi|^{1/\sigma})
\end{equation}
for all $\xi\in\R^n$, $|\xi|\geq1$. To verify that
$u\in W^{2,1}([0,T],\gamma^{(\sigma)})$ it is sufficient to show
that there exist
$M',\delta'>0$ such that
\begin{equation}\label{eq2.5}
|v(t,\xi)|^2+|v'(t,\xi)|^2\leq M'\exp(-\delta'|\xi|^{1/ \sigma})
\end{equation}
for all $t\in[0,T]$ and for all $\xi\in\R^n$, $|\xi|\geq1$. We consider
first the case $t_0=T$. For $\eps\in]0,T]$, we set
\begin{equation}\label{aeps1}
a_\eps(t,\xi):=
\begin{cases}
a(t,\xi)+\eps^{2-(r+1/s)}(T-t)^{-2}&\text{for $0\leq t\leq T-\eps$}\\
\eps^{-(z+r)}a(t,\xi)(T-t)^{z+r}+\eps^{-(r+1/s)}&\text{for $T-\eps
\leq t\leq T$}
\end{cases}
\end{equation}
where $z$ is any positive number such that
\begin{equation}\label{condz1}
z>\max\left\{1/s, \left(p+1/q\right)-r-1\right\}.
\end{equation}
Then
\begin{equation}\label{x1}
a_\eps(t,\xi)=
\begin{cases}
\alpha(t,\xi)(T-t)^{-r}+\eps^{2-(r+1/s)}(T-t)^{-2}
&\text{for $0\leq t\leq T-\eps$}\\
\eps^{-(z+r)}\alpha(t,\xi)(T-t)^{z}+\eps^{-(r+1/s)}
&\text{for $T-\eps\leq t\leq T$}
\end{cases}
\end{equation}
and
\begin{equation}\label{x2}
a'_\eps(t,\xi)
=\begin{cases}
\beta(t,\xi)(T-t)^{-p}-2\eps^{2-(r+1/s)}(T-t)^{-3}&
\text{for $0\leq t\leq T-\eps$}\\
\eps^{-(z+r)}\left(\beta(t,\xi)(T-t)^{z+r-p}\right.&\\
\quad\quad\quad\quad\quad\,\,\left.-(z+r)\alpha(t,\xi)(T-t)^{z-1}\right)
&\text{for $T-\eps\leq t\leq T$}\\
\end{cases}
\end{equation}
Our choice of $z$ implies that $a_{\epsilon}(\cdot,\xi)\in W^{1,1}(0,T)$
for all
$\xi\in\R^n\setminus\{0\}$. By (\ref{x1}) and (\ref{x2}) we get
\begin{align*}
\int_0^T\frac{|a'_\eps(t,\xi|}{a_\eps(t,\xi)}\,d t
&\leq\int_0^{T-\eps}\frac{|\beta(t,\xi)|
(T-t)^{-p}}{\eps^{2-(r+1/s)}(T-t)^{-2}}\,d t\\
&+\int_0^{T-\eps}\frac{
2\eps^{2-(r+1/s)}(T-t)^{-3}}{\eps^{2-(r+1/s)}(T-t)^{-2}}\,d t\\
&+\int_{T-\eps}^T\frac{\eps^{-(z+r)}|\beta(t,\xi)|(T-t)^{z+r-p}}
{\eps^{-(r+1/s)}}\,d t\\
&+\int_{T-\eps}^T\frac{\eps^{-(z+r)}
(z+r)|\alpha(t,\xi)|(T-t)^{z-1}}{\eps^{-(r+1/s)}}\,d t
\end{align*}
The choice of $z$ allows us to use H\"older inequality; an easy
computation shows that
\begin{equation}\label{est1}
\int_0^T\frac{|a'_\eps(t,\xi|}{a_\eps(t,\xi)}\,d t
\leq C'(1+|\log\eps|)\eps^{-(p+1/q)+(r+1/s)+1},
\end{equation}
where $C'$ is a constant depending only on $C,r,s,p,q$ and $z$.
On the other hand,
\begin{align*}
\int_0^T\frac{|a_\eps(t,\xi)-a(t,\xi)|}{a_\eps(t,\xi)^{1/2}}\,d t
&=\int_0^{T-\eps} \frac {\eps^{2-(r+1/s)}(T-t)^{-2}}
{\eps^{1-(1/2)(r+1/s)}(T-t)^{-1}} \, dt \\
&+\int_{T-\eps}^{T}\frac {\eps^{-(z+r)}\alpha(t,\xi)(T-t)^z}
{\eps^{-(1/2)(r+1/s)}}  \, dt\\
&+\int_{T-\eps}^{T}\frac {\eps^{-(r+1/s)}}{\eps^{-(1/2)(r+1/s)}}  \, dt
+\int_{T-\eps}^{T}\frac {\alpha(t,\xi)(T-t)^{-r}}{\eps^{-(1/2)(r+1/s)}}  \, dt.
\end{align*}
The first three summands on the right hand side can be estimated again
by using H\"older inequality. In order to estimate the fourth summand, we
shall
distinguish the case $(r+1/s)<1$ and $(r+1/s)=1$. In the first case, we use
once more H\"older inequality; in the second case, we use the fact that
$\alpha(t,\xi)(T-t)^{-r}=a(t,\xi)\in L^1(0,T)$. At the end, we get
\begin{equation}\label{est2}
\int_0^T\frac{|a_\eps(t,\xi)-a(t,\xi)|}{a_\eps(t,\xi)^{1/2}}\,d t
\leq C''(1+|\log\eps|)\eps^{-(1/2)(r+1/s)+1},
\end{equation}
where $C''$ is a constant depending only on $C,r,s,p,q$ and $z$.
By (\ref{eq2.3}), (\ref{est1}) and (\ref{est2}) we obtain
\begin{multline}\label{energyest1}
E(t,\xi)\\
\leq E(0,\xi)\exp{(\tilde C(1+|\log\eps|)(\eps^{-(p+1/q)+(r+1/s)+1}
+|\xi|\eps^{-(1/2)(r+1/s)+1}))}
\end{multline}
for all $t\in[0,T]$ and for all $\xi\in\R^n$, $\xi\geq 1$, where
$\tilde C$ is a positive constant depending only on $C,r,s,p,q$ and $z$.

Now, by (\ref{eq2.2}) and (\ref{aeps1}), we have
\begin{equation}
E_\eps(0,\xi)\leq\left(a(0,\xi)+T^{-(r+1/s)}\right)|\xi|^2|v(0,\xi)|^2
+|v'(0,\xi)|^2
\end{equation}
and
\begin{equation}
E_\eps(t,\xi)\geq T^{-2}\eps^{2-(r+1/s)}|\xi|^2|v(t,\xi)|^2
+|v'(t,\xi)|^2.
\end{equation}
Then choosing $\eps:=|\xi|^{-[{(p+1/q)-\frac32(r+1/s)}]^{-1}}$ we deduce
\begin{multline*}
T^{-2}|\xi|^{2-\frac{2-(r+1/s)}{(p+1/q)-\frac32(r+1/s)}}|v(t,\xi)|^2+
|v'(t,\xi)|^2\\
\leq (\tilde K|\xi|^2|v(0,\xi)|^2
+|v'(0,\xi)|^2)
\exp{(\tilde C(1+|\log|\xi|)|\xi|^{\frac{(p+1/q)-(r+1/s)-1}
{(p+1/q)-\frac32(r+1/s)}})}.
\end{multline*}
Using the Paley-Wiener theorem, the well-posedness follows for all
$1\leq\sigma<
\frac{(p+1/q)-\frac32(r+1/s)}{(p+1/q)-(r+1/s)-1}$.

\par If $t_0=0$, for $\eps\in]0,T]$ we set
\begin{equation}\label{aeps2}
a_\eps(t,\xi):=
\begin{cases}
\eps^{-(z+r)}a(t,\xi)t^{z+r}+\eps^{-(r+1/s)}&\text{for $0
\leq t\leq \eps$}\\
a(t,\xi)+\eps^{2-(r+1/s)}t^{-2}&\text{for $\eps\leq t\leq T$}\\
\end{cases}
\end{equation}
where $z$ satisfies (\ref{condz1}). Our choice of $z$ implies that
$a_{\epsilon}(\cdot,\xi)\in W^{1,1}(0,T)$ for all $\xi\in\R^n\setminus\{0\}$.
So, in particular, $a_{\epsilon}(\cdot,\xi)$ is continuous on $[0,T]$.
Arguing as before, we obtain (\ref{energyest1}). An easy computation
shows that $|a(t,\xi)|\leq \tilde K t^{1-(p+1/q)}$ for all
$\xi\in\R^n\setminus\{0\}$. It follows that
\begin{multline*}
a_\eps(0,\xi)=\lim_{\tau\to 0}a_\eps(\tau,\xi)
=\lim_{\tau\to 0}(\eps^{-(z+r)}a(\tau,\xi)\tau^{z+r}+\eps^{-(r+1/s)})\\
\leq\tilde K \limsup_{\tau\to 0}
(\eps^{-(z+r)}\tau^{z+r+1-(p+1/q)}+\eps^{-(r+1/s)}).
\end{multline*}
By (\ref{condz1}) we deduce that $a_\eps(0,\xi)\leq \tilde K \eps^{-(r+1/s)}$.
It follows that
\begin{equation}
E_\eps(0,\xi)\leq \tilde K\eps^{-(r+1/s)}|\xi|^2|v(0,\xi)|^2+|v'(0,\xi|^2.
\end{equation}
Moreover, we have also
\begin{equation}
E_\eps(t,\xi)\geq T^{-2}\eps^{2-(r+1/s)}|\xi|^2|v(t,\xi)|^2+|v'(t,\xi)|^2.
\end{equation}
Then, choosing again $\eps:=|\xi|^{-[{(p+1/q)-\frac32(r+1/s)}]^{-1}}$,
we deduce
\begin{multline*}
|\xi|^{2-\frac{2-(r+1/s)}{(p+1/q)-\frac32(r+1/s)}}|v(t,\xi)|^2+
|v'(t,\xi)|^2
\leq (\tilde K|\xi|^{2+\frac{(r+1/s)}{(p+1/q)-\frac32(r+1/s)}}|v(0,\xi)|^2\\
+|v'(0,\xi)|^2)
\exp{(\tilde C(1+|\log|\xi|)|\xi|^{\frac{(p+1/q)-(r+1/s)-1}
{(p+1/q)-\frac32(r+1/s)}})}.
\end{multline*}
Using the Paley-Wiener theorem, the well-posedness follows again for all
$1\leq\sigma<
\frac{(p+1/q)-\frac32(r+1/s)}{(p+1/q)-(r+1/s)-1}$.
\par Finally, if $t_0\in]0,T[$, it will be sufficient to solve first
the Cauchy problem in $[0,t_0]$, then to solve the problem
in $[t_0,T]$ with the initial data obtained from the previous one
and finally to glue together the two solutions. \end{proof}

\par In order to prove Theorem \ref{th2}, we proceed exactly like in the
proof of Theorem \ref{th1}. In this case the role of
condition (\ref{supplcond}) is played by the estimate
\begin{equation}
a(t,\xi)\leq C'|t-t_0|^{-(p+1/q)+1}\quad\text{for all
$\xi\in\R^n\setminus\{0\}$,}
\end{equation}
which is a direct consequence of condition (H2).
The function $a_\eps(\cdot,\xi)$ is defined
by
\begin{equation}\label{aeps3}
a_\eps(t,\xi):=
\begin{cases}
a(t,\xi)+\eps^{3-(p+1/q)}(T-t)^{-2}&\text{for $0\leq t\leq T-\eps$}\\
a(T-\eps,\xi)+\eps^{1-(p+1/q)}&\text{for $T-\eps\leq t\leq T$}
\end{cases}
\end{equation}
if $t_0=T$ and by
\begin{equation}\label{aeps4}
a_\eps(t,\xi):=
\begin{cases}
a(\eps,\xi)+\eps^{1-(p+1/q)}&\text{for $0
\leq t\leq \eps$}\\
a(t,\xi)+\eps^{3-(p+1/q)}t^{-2}&\text{for $\eps\leq t\leq T$}\\
\end{cases}
\end{equation}
if $t_0=0$. Arguing like in the proof of Theorem \ref{th1}, we get
\begin{equation}
\int_0^T\frac{|a'_\eps(t,\xi|}{a_\eps(t,\xi)}\,d t
\leq C''(1+|\log\eps|)\eps^{(p+1/q)-3},
\end{equation}
and
\begin{equation}
\int_0^T\frac{|a_\eps(t,\xi)-a(t,\xi)|}{a_\eps(t,\xi)^{1/2}}\,d t
\leq C''(1+|\log\eps|)\eps^{-(1/2)(p+1/q)+3/2}
\end{equation}
and the conclusion follows by choosing $\eps:=|\xi|^{-(2/3)[3-(p+1/q)]^{-1}}$.

\par Theorem \ref{th3} is the strictly hyperbolic version of
Theorem \ref{th2}. We define again $a_\eps$ by (\ref{aeps3}) and (\ref{aeps4}),
but in this case the positive lower bound for $a(t,\xi)$ allows us to obtain
better estimates  for $\int_0^T\frac{|a'_\eps(t,\xi|}{a_\eps(t,\xi)}\,d t$.
Let us consider, for example, the case $t_0=T$.
First  observe that, by rescaling the $x$ variable if
necessary, we can always assume that $\lambda_0=1$. Then we can minorize
$a_\eps(t,\xi)$ by the constant $1$ on $[0,T-\eps^{(1/2)[3-(p+1/q)]}]$ and by
$\eps^{3-(p+1/q)}(T-t)^{-2}$ on $[T-\eps^{(1/2)[3-(p+1/q)]},T-\eps]$. So we
obtain
that
\begin{equation}
\int_0^T\frac{|a'_\eps(t,\xi|}{a_\eps(t,\xi)}\,d t
\leq C''(1+|\log\eps|)\eps^{(1/2)((p+1/q)-1)((p+1/q)-3)}.
\end{equation}
The conclusion follows by choosing
$\eps:=|\xi|^{-2[p+1/q]^{-1}[3-(p+1/q)]^{-1}}$.

\par Finally, we give the
\begin{proof}[Proof of Theorem \ref{th4}]
Since $u_0$, $u_1\in C^\infty_0$, the Paley-Wiener
theorem ensures that for all $\zeta>0$ there exists $M_\zeta>0$ such that
\begin{equation}\label{cinf1}
|v(0,\xi)|^2+|v'(0,\xi)|^2\leq M_\zeta |\xi|^{-\zeta}
\end{equation}
for all $\xi\in\R^n$, $|\xi|\geq1$. To verify that
$u\in W^{2,1}([0,T],C^\infty_0)$ it is sufficient to show  that for all
$\eta>0$
there exists
$M_\eta>0$ such that
\begin{equation}\label{cinf2}
|v(t,\xi)|^2+|v'(t,\xi)|^2\leq M_\eta |\xi|^{-\eta}
\end{equation}
for all $t\in[0,T]$ and for all $\xi\in\R^n$, $|\xi|\geq1$.
We give the details only in the case $t_0=T$.
If $q=1$, then necessarily $p=0$. This means that
$a(\cdot,\xi)\in W^{1,1}(0,T)$ and it is well known that this is enough to
detect $C^\infty$-well-posedness of the Cauchy problem
(\ref{eq1}), (\ref{eq2}). If $q>1$, for $\eps\in]0,T]$, we set
\begin{equation}\label{aeps5}
a_\eps(t,\xi):=
\begin{cases}
a(t,\xi)&\text{for $0\leq t\leq T-\eps$}\\
a(T-\eps,\xi)&\text{for $T-\eps\leq t\leq T$}
\end{cases}
\end{equation}
Now observe that
\begin{multline*}
|a(t,\xi)|\leq|a(0,\xi)|+\int_0^t|a'(\tau,\xi)|\,d\tau\leq
|a(0,\xi)|+\int_0^t\beta(\tau,\xi)(T-\tau)^{-p}\,d\tau\\
\leq|a(0,\xi)|+
\|\beta(\cdot,\xi)\|_{L^q}\left(\int_0^t(T-\tau)^{-pq'}\,d\tau\right)^{1/q'}\leq
C(1+|\log(T-t)|^{1/q'})
\end{multline*}
An easy computation shows that
\begin{equation}
\int_0^T|a'_\eps(t,\xi)|\,dt\leq C'|\log\eps|^{1/q'}
\end{equation}
and
\begin{equation}
\int_0^T|a_\eps(t,\xi)-a(t,\xi)|\,dt\leq C'\eps|\log\eps|^{1/q'}.
\end{equation}
Then we deduce by (\ref{eq2.3}) that
\begin{multline}
|\xi|^2|v(t,\xi)|^2+|v'(t,\xi)|^2\\
\leq(a(0,\xi)|\xi|^2|v(0,\xi)|^2+v'(0,\xi)|^2)
\exp(C'|\log\eps|^{1/q'}+C'|\xi|\eps|\log\eps|^{1/q'}).
\end{multline}
Here, for simplicity, we have assumed that $\lambda_0=1$.
Choosing $\eps:=|\xi|^{-1}$, we obtain
\begin{multline}
|\xi|^2|v(t,\xi)|^2+|v'(t,\xi)|^2\\
\leq(a(0,\xi)|\xi|^2|v(0,\xi)|^2+v'(0,\xi)|^2)
\exp(C'|\log|\xi||^{1/q'}).
\end{multline}
Now, for $|\xi|\geq e$, we have $|\log|\xi||^{1/q'}\leq|\log|\xi||$, and
hence
\begin{equation}
|\xi|^2|v(t,\xi)|^2+|v'(t,\xi)|^2
\leq (a(0,\xi)|\xi|^2|v(0,\xi)|^2+v'(0,\xi)|^2) |\xi|^{C'}.
\end{equation}
By the Paley-Wiener theorem, the well-posedness in $C^\infty_0$ follows.
\end{proof}

\end{document}